\documentstyle[12pt]{article}

\begin{document}
\newcommand{\lra}{\longrightarrow}
\newcommand{\ra}{\rightarrow}
\newcommand{\ga}{\gamma}
\newcommand{\st}{\stackrel}
\newcommand{\ol}{\overline}
\newcommand{\ul}{\underline}
\newcommand{\equ}{\equiv}

\title{ \Large\bf Subgroup Theorems for the Baer-invariant of Groups }
\author{M.R.R.M{\small\bf OGHADDAM}, B.M{\small\bf ASHAYEKHY, AND}
S.K{\small\bf AYVANFAR}\\Department
of Mathematics, Mashhad Univesity\\P.O.Box 1159-91775, Mashhad, Iran}
\date{}

\maketitle

\begin{abstract}
 M.R.Jones and J.Wiegold in [3] have shown that if $G$ is a finite group with a
subgroup $H$ of finite index $n$ , then the $n$-th power of Schur multiplier
of $G$ , $M(G)^n$ , is isomorphic to a subgroup of $M(H)$ .

 In this paper we prove a similar result for the centre by centre by $w$
variety of groups, where $w$ is any outer commutator word. Then
using a result of M.R.R.Moghaddam [6], we will be able to deduce
a result of Schur's type ( see [4,9] ) with respect to the
variety of nilpotent groups of class at most $c$ $(c\geq 1)$,
when $c+1$ is any prime number or $4$.\\
A.M.S.Classification(1991): 20E07,20E10,20F12
\end{abstract}
\ \ \\
 \begin{center}  {\bf 1. I{\small NTRODUCTION}}\end{center}

 Let $F_{\infty}$ be the free group freely generated by an infinite set
$\{x_1,x_2,\ldots\}$ . Let $w$ be an outer commutator word in $F_{\infty}$ ,
and assume ${\cal V}$ to be the variety of groups defined by the commutator
word $v=[w,x_1,x_2]$ . Let $G$ be a finite  group with a ${\cal V}$- stem
cover $G^*/L\cong G$,say. If $H$ is a subgroup of $G$ of finite index $n$ such
that $H\cong  B/L$, then it is shown (Theorem 3.6) that the commutator
subgroup $[w(G), {\cal V}M(G)^n]$ is isomorphic to a subgroup of ${\cal V}M(H)$
, where $w(G)$ is the verbal subgroup of $G$ with respect to the word $w$ ,
and ${\cal V}M(G)$ is the Baer-invariant of $G$ with respect to the variety
of groups $\cal V$ .

 Now, if we assume $w$ to be an empty word, then we obtain the result of
M.R.Jones and J.Wiegold [3] . To prove the above fact, one needs to have the
concept of the transfer of a group (see [8]) in the above variety $\cal V$ ,
which has been done in the Definition 3.3, (see also Proposition 3.4).

 I.Schur in 1904 [9] showed that if $G$ is a finite group and $G_p$ is a Sylow
$p$-subgroup of $G$, then the Sylow $p$-subgroup of the Schur multiplier of
$G$, denoted by $M(G)_p$ , is
\\ \_\hrulefill\_\ \ \ \ \ \ \ \
\ \ \ \ \ \ \ \ \ \ \ \ \ \ \ \ \ \ \ \ \ \ \ \ \ \ \ \ \ \ \ \ \ \ \ \ \ \  \\
Communication address:c/ M.R.R.Moghaddm, Department of Mathematics, Ferdowsi
University
of Mashhad, P.O.Box 1159-91775, Mashhad, Iran. Fax:05-817749 Tel:0098-51-
817749, e-mail:rezarm@toos.um.ac.ir\ \  or\ \ rezarm@irearn.bitnet\\
isomorphic to a subgroup of
the Schur multiplier of its Sylow $p$-subgroup, $M(G_p)$ .

 In section 4, by using our previous result together with a result of
M.R.R.Moghaddam (see Theorem 2.3 [6]), it is proved a theorem of Schur's type
(Theorem 4.3) with respect to the variety of nilpotent groups, ${\cal N}_c$ ,
of class at most $c$ , where $c+1$ is any prime number or $c=3$ .

 Finally, in section 5 we give an example to show that the main result can not
be futher generalized. In other words, it is in best possible shape. \\

\begin{center}{\bf 2. N{\small\bf OTATION AND} P{\small\bf RELIMINARIES}}
\end{center}

 Let $x,y$ be two elements of a group $G$ then $[x,y]$ , the {\it commutator}
of $x$ and
$y$ , and $x^y$ denote the elements $x^{-1}y^{-1}xy$ and $y^{-1}xy$
, respectively.
The commutator of higher weight is defined inductively as follows:
$$[x_1,x_2,\ldots,x_n]=[[x_1,x_2,\ldots,x_{n-1}],x_n]\ .$$

 If $H$ and $K$ are two subgroups of a group $G$ then $[H,K]$ denotes the
subgroup of
$G$ generated by all the commutators $[h,k]$ with $h$ in $H$ and $k$ in $K$ .
In particular,
if $H=K=G$ then $[G,G]$ , which is denoted by $G'$ , is the {\it derived
subgroup} of $G$ . The lower and upper centeral series are denoted by
$\gamma_c(G)$\ \  $(c\geq1)$ and $Z_d(G)$\ \  $(d\geq 0)$ , respectively.

 Let $F_{\infty}$ be the free group freely generated by an infinite countable
set $\{x_1,x_2,\ldots \}$. If
$u=u(x_1,\ldots ,x_s)$ and $v=v(x_1,\ldots ,x_s)$ are two words in $F_{\infty}$
, then the {\it composite} of $u$ and $v$ , $u\circ v$ , is defined as follows:
$$u\circ v=u(v(x_1,\ldots ,x_t),\ldots,v(x_{(s-1)t+1},\ldots,x_{st}))\ \ .$$

 In particular, the composite of some nilpotent words is called {\it
polynilpotent word}, i.e. $$\gamma_{c_1+1,\ldots,c_t+1}=\gamma_{c_1+1}\circ
\gamma_{c_2+1}\circ \cdots \circ \gamma_{c_t+1}\ \ ,$$
where $\gamma_{c_i+1}$\ \  $(1\leq i\leq t)$ is a nilpotent word in distinct
variables.

{\it Outer commutator words} are defined inductively, as follows:\\
 The word $x_i$ is an outer commutator word (henceforth o.c.word) of weight
one. If $u=u(x_1.\ldots,x_s)$ and $v=v(x_{s+1},\ldots,x_{s+t})$ are o.c.words
of weights $s$ and $t$ , respectively, then
$$w(x_1,\ldots,x_{s+t})=[u(x_1,\ldots,x_s),
v(x_{s+1},\ldots,x_{s+t})]$$ is an o.c.word of weight $s+t$ .

 Let $V$ be a subset of $F_{\infty}$ and ${\cal V}$ be the variety of groups
defined by the set of laws $V$ . Then following P.Hall [2] in 1957 , introduce
the notion of the {\it verbal subgroup}, $V(G)$ , and the {\it marginal
subgroup}, $V^*(G)$ , associated with the variety ${\cal V}$ and the group $G$,
which we assume that the reader is familiar with these concepts.

 If $G$ is a group with a free presentation
$ 1 \lra R \lra F \lra G \lra 1\ \ ,$ then the {\it
Baer\_invariant} of $G$ , with respect to the variety ${\cal V}$ , is defined
to be $${\cal V}M(G)=\frac {R\cap V(F)} {[RV^*F]}\ \ \ ,$$
where $V(F)$ is the verbal subgroup of $F$ and
$$[RV^*F]=<v(f_1,\ldots,f_ir,\ldots,
f_s)v(f_1,\ldots,f_s)^{-1}\mid v\in V,f_i\in F,r\in R,1\leq i\leq s>.$$

One notes that the Baer-invariant of the group $G$ is always abelian and
independent of the choice of the free presentation of $G$ (see
M.R.R.Moghaddam [6] or C.R.Leedham-Green and S.McKay [5]).

In particular, if ${\cal V}$ is the variety of abelian or nilpotent groups of
class at most $c$ $(c\geq 2)$ , then the Baer-invariant of $G$ will be $R\cap
F'/ [R,F]$ , which is isomorphic to the Schur-multiplier of $G$
(I.Schur [9]), or $R\cap \gamma_{c+1}(F) / [R,_cF] $, respectively
(M.R.R.Moghaddam [6]).

Let ${\cal V}$ be a variety of groups and $G$ be an arbitrary group, then an
extension $$ 1 \lra L \lra E \lra G \lra 1 \ \ \ \ \ \ \ (*)$$
is called a ${\cal V}$-{\it stem extension} of $G$ , if $L\subseteq V(E)\cap
V^*(E)\ .$ The extension $(*)$ will be called ${\cal V}$-{\it stem cover} of
$G$ , when we also have $L\cong {\cal V}M(G)$ , and in this case $E$ is said to
be a ${\cal V}$-{\it covering group} of $G$ (see [5]).\\
{\bf Definition 2.1}

 Let ${\cal V}$ be a variety of groups defined by the set of laws $V$ . Then
${\cal V}$ is said to be a {\it Schur-Baer variety}; whenever, the marginal
factor group $G/V^*(G)$ is finite of order $m$ , say, then the verbal subgroup
$V(G)$ is also finite and $|V(G)|\ |\ m^k \ \ ,\ k\in {\bf N}$ , for arbitrary
group $G$ .

  Among many other results, I.Schur [9] showed that the variety of abelian
groups has Schur-Baer property. Also R.Baer [1] proved that the variety
defined by o.c. words has the same property.

The following results are needed later. \\
{\bf Lemma 2.2}

 Let ${\cal V}$ be a variety defined by the set of laws $V$ , then the
following conditions are equvalent:

$(i)$ ${\cal V}$ is a Schur-Baer variety.

$(ii)$ For any finite group $G$ , the Baer-invariant ${\cal V}$M(G) is of
order dividing a power of $|G|$ .\\
{\bf Proof.}  See [5] . \\
{\bf Theorem 2.3}  (M.R.R.Moghaddam [6])

 Let ${\cal N}_c$ be the variety of nilpotent groups of class at most $c$ ,and
let $G_1$ and $G_2$ be two arbitrary groups. Then

$(i)$ ${\cal N}_cM(G_1\times G_2)\cong {\cal N}_cM(G_1)\oplus {\cal
N}_cM(G_2)\oplus T(G_2,G_1)_{c+1}$, when $c+1$ is a prime number.

$(ii)$ ${\cal N}_3M(G_1\times G_2)\cong {\cal N}_3M(G_1)\oplus {\cal
N}_3M(G_2)\oplus \hat {T}(G_2,G_1)_4$\ \ ,\\ where $T(G_2,G_1)_{c+1}$ denotes
the summation of all tensor product corresponding to the basic commutators of
weight $c+1$ on two letters, and
$$\hat {T}(G_2,G_1)_4=(\overline {G_2}\otimes \overline {G_1} \otimes \overline
{G_1}\otimes \overline {G_1} )\ \oplus \  \sum _{i=1}^2(\overline {G_2}\otimes
\overline {G_1}\otimes \overline {G_i}\otimes \overline {G_1}) \\ \oplus \
(\overline
{G_2}\otimes \overline {G_1})\wedge (\overline {G_2}\otimes \overline {G_1})\
\ ,$$
where $\overline {G}_i=(G_i)_{ab}\ ,\ i=1,2$ , and $\wedge$ is the wedge
product. \\
{\bf Proof.}  See [6] . \\

\begin{center}{\bf 3. T{\small\bf HE} M{\small\bf AIN} R{\small\bf ESULTS}}
\end{center}

 In 1973 , M.R.Jones and J.Wiegold [3] showed that if $G$ is a finite
group, and $H$ a subgroup of $G$ of index $n$ , then the $n$-th power of
$M(G)$ , i.e. $M(G)^n$, is
isomorphic to a subgroup of $M(H)$ , which generalizes Schur's Theorem [4,9].

 In this section the above result will be generalized to the centre by centre
by $w$ variety of groups, where $w$ is any outer commutator word.

 We assume $G$ to be a finite group, unless otherwise stated.\\
{\bf Lemma 3.1}

 Let ${\cal V}$ be a variety of groups, and $H$ be a group with a ${\cal
V}$-marginal subgroup $A$ , i.e. $A\leq V^{*}(H)\ .$ If $G=H/A$ , then
$V(H)\cap A$ is a homomorphic image of ${\cal V}M(G)\ .$\\
{\bf Proof.}

 Let $H\cong F/R$ be a free presentation of $H$ , then for some subgroup $T$ of
$F$ , we have $A\cong T/R$ , hence $G\cong F/T$ , and by the definition
$${\cal V}M(G)=\frac {T\cap V(F)} {[TV^{*}F]}\ \ ,$$
and using Dedekind's law,
$$V(H)\cap A=\frac {V(F)R} {R}\cap \frac {T} {R} =\frac {V(F)R\cap T} {R} =
\frac {(V(F)\cap T)R} {R}$$ $$\cong \frac {V(F)\cap T} {(V(F)\cap T)\cap R} =
\frac {V(F)\cap T} {V(F)\cap R}\ \ .$$
Since $T/R$ is a marginal subgroup of $F/R$ , it follows that
$[TV^{*}F]\subseteq R$ . So we have
$$V(H)\cap A=\frac {V(F)\cap T} {V(F)\cap R}\cong \frac {V(F)\cap T/
 [TV^{*}F]}{V(F)\cap R/ [TV^{*}F]}\ \ .$$
Thus the result holds.\ \ $\Box$ \\
{\bf Corollary 3.2}

 Let ${\cal V}$ be a Schur-Baer variety and $G$ be a finite group. If $H$ is
any  other group with a ${\cal V}$-marginal subgroup $A$ such that $H/A\cong
G$ , then $V(H)\cap A$ is isomorphic to a subgroup of ${\cal V}M(G)\ .$\\
{\bf Proof.}

 Using Lemma 3.1 , $V(H)\cap A$ is a homomorphic image
of  ${\cal V}M(G)\ .$ By Lemma 2.2, since ${\cal V}$ is a Schur-Baer variety
and $G$ is finite, it implies that ${\cal V}M(G)$ is a finite abelian group.
Now, by the property of finite abelian groups, $V(H)\cap A$ is isomorphic to a
subgroup of ${\cal V}M(G)\ .$ \ \ $\Box$

 Let ${\cal W}$ be a variety defined by an outer commutator word
$w=w(x_1,\ldots ,x_r)$ , and ${\cal V}$ be the centre by centre by
${\cal W}$ variety of groups. Clearly the variety ${\cal V}$ will be
defined by the commutator word   $[w,x_{1},x_{2}]\ .$ Throughout the rest of
this section we always work with the above variety ${\cal V}\ .$

 The following definition and proposition are very important in our
investigation, in which the concept of the {\it tranfer } has been generalized
to the above variety ${\cal V}\ .$\\
{\bf Definition 3.3}

 Let $E$ be any group with a subgroup $B$ of finite index $n$ , and
$\{t_1,\ldots ,t_n\}$ be a right transversal of $B$ in $E$ . Let ${\cal V}$ be
the variety of groups as defined before, and
$$\pi :B\longrightarrow B/V(B)$$  $$\  \ b\longmapsto \overline {b}\ \ \ \ \
\ \ $$ be the natural homomorphism, where $V(B)$ is the verbal subgroup of $B$
with respect to the variety ${\cal V}\ .$

 Now, for all $y_1,\ldots,y_r$ in $B$ and $x$ in $E$ we define the following
map $$\tau _{y_1,\ldots ,y_r}=\tau :E\longrightarrow B/V(B)\ \ \ \ \ \ , $$
given by \ \ \ $x\longmapsto \prod_{i=1}^n [w(\overline {y_1},\ldots ,\overline
{y_r}),(t_ixt_{ix}^{-1})\pi ]\ ,$\ \
where $Bt_ix=Bt_{ix}$ , for all $x\in E\ .$

 The following result shows that the map $\tau $ has the property of being a
generalized transfer of the group $E$\ .\\
{\bf Proposition 3.4}

 With the above set up, $\tau $ is a well-defined homomorphism and indepedent
of the choice of the transversal set. Also, if $x$ is marginal in $E$, then
$$(x)\tau =\overline {[w(y_1,\ldots ,y_r),x^n]}\ \ .$$
{\bf Proof.}

 Clearly $\tau $ is well-defined. Let $x,y\in E$, then
we have \ \ \  $t_ixyt_{ixy}^{-1}=t_ixt_{ix}^{-1}t_{ix}yt_{(ix)y}^{-1}\ .$\\
By the form of the commutator word which defines the variety ${\cal V}$, modulo
$V(B)$ , we obtain
$$[w(a_1,\ldots ,a_r),[a_{r+1},a_{r+2}]]\equiv [w(a_1,\ldots
,a_r),a_{r+1},a_{r+2}]\equiv 1\ \ ,\ \ (*) $$   for all $a_i\in B$ . Thus
$$(xy)\tau =\prod_{i=1}^n[w(\overline {y_1},\ldots ,\overline
{y_r}),(t_ixyt_{ixy}^{-1})\pi]$$
$$\     \ =\prod_{i=1}^n[w(\overline {y_1},\ldots ,\overline
{y_r}),(t_ixt_{ix}^{-1})\pi (t_{ix}yt_{(ix)y}^{-1})\pi ]$$
$$\ \ \ \ \     \ \equiv \prod_{i=1}^n[w(\overline {y_1},\ldots ,\overline
{y_r}),(t_ixt_{ix}^{-1})\pi ] [w(\overline {y_1},\ldots ,\overline
{y_r}),(t_jyt_{jy}^{-1})\pi ] \ \ \ \ \ \ \  (mod\ \ V(B))$$
$$=(x)\tau (y)\tau \ \ \ \ , \ \ \ {\it putting}\ \ ix=y\ \ .$$
Hence $\tau $ is a homomorphism.

 Now, let $\{t_1^{\prime},\ldots ,t_n^{\prime}\}$ be another transversal set of
$B$ in $E$, then $t_i^{\prime}=b_it_i$ for some $b_i\in B\ ,\ \ 1\leq i\leq n\
.$ This implies that, for all $x\in E\ .$
$$t_i^{\prime}x{t_{ix}^{\prime}}^{-1}=b_it_ixt_{ix}^{-1}b_{ix}^{-1}\equiv
b_ib_{ix}^{-1}t_ixt_{ix}^{-1}\ \ \ \ (mod\ \ \gamma_2(B))\ .$$
Now, using $(*)$ and the above fact, mod $V(B)$ ,we have
$$\prod_{i=1}^n[w(\overline {y_1},\ldots ,\overline
{y_r}),(t_i^{\prime}x{t_{ix}^{\prime}}^{-1})\pi]\equiv [w(\overline
{y_1},\ldots ,\overline
{y_r}),\prod_{i=1}^n(t_i^{\prime}x{t_{ix}^{\prime}}^{-1})\pi ]$$
$$\ \ \ \ \ \ \ \ \ \ \ \equiv [w(\overline {y_1},\ldots ,\overline
{y_r}),\prod_{i=1}^{n} (t_ixt_{ix}^{-1})\pi]\ \ \ \ \ \ \ ({\it since}\ \
\prod_{i=1}^nb_ib_{ix}^{-1}=1)$$
$$ \equiv \prod_{i=1}^n[w(\overline {y_1},\ldots ,\overline
{y_r}),(t_ixt_{ix}^{-1})\pi]=(x)\tau ,\ \ \ \ \ \ \ \ \ \ \ \ $$
which implies that $\tau $ is independent of the choice of the transversal
set.

 Clearly every elements $x$ of $E$ acts on the transversal set $\{t_1,\ldots
, t_n\}$ . Assume there are $k$ orbits with the following action:
$$t_{i_{1}}\ ,\ t_{i_{2}}=t_{i_{1}}x\ ,\ \ldots
,\ t_{i_{l_{i}}}=t_{i_{1}}x^{l_i-1}\ ,\ \
{\it and}\ \ t_{i_{1}}=t_{i_{1}}x^{l_i}\ ,$$  $ where\ \ 1\leq i\leq k\ \
,l_1+\ldots +l_k=n\ .$ Clearly
$$t_{i_{j}}x=\left \{ \begin{array}{lc} t_{i_{j+1}} & ,\ if\ \ j\not =l_i\\
 t_{i_{1}} & ,\ if\ \ j=l_i\end{array} \right.$$
which implies that
$$t_{i_{j}}xt_{i_{j}x}^{-1}=\left \{ \begin{array}{lc} 1 & ,\ if\ \ j\not
=l_i\\ t_{i_{1}}x^{l_i}t_{i_{1}}^{-1} & ,\ if\ \ j=l_i\end{array} \right.$$
Hence
$$(x)\tau =\prod_{i=1}^n[w(\overline {y_1},\ldots ,\overline
{y_r}),(t_ixt_{ix}^{-1})\pi ]$$
$$\ \ \ \ \ =\prod_{i=1}^k[w(\overline {y_1},\ldots ,\overline
{y_r}),(t_{i_{1}}x^{l_i}t_{i_{1}}^{-1})\pi ]\ \ .$$

 Now, if $x$ is marginal in $E$ , using Turner-Smith's Theorem [10], and
Witt's identity we have
$$ 1=[w(x_1,\ldots ,x_r),x_{r+1},x]=[w(x_1,\ldots
,x_r),x,x_{r+1}]$$ $$=[w(x_1,\ldots
,x_r),[x,x_{r+1}]]\ \ \ ,\ {\it for\ all} \ \ x_i\in E\ \ .$$
So the following holds:
$$ (x)\tau =\prod_{i=1}^k[w(\overline {y_1},\ldots ,\overline
{y_r}),(t_{i_{1}}x^{l_i}t_{i_{1}}^{-1})\pi ]$$
$$ =\overline {\prod_{i=1}^k[w(y_1,\ldots
,y_r),t_{i_{1}}t_{i_{1}}^{-1}x^{l_i}[x^{l_i},t_{i_{1}}^{-1}]]}$$
$$ =\overline {\prod_{i=1}^k[w(y_1,\ldots
,y_r),[x^{l_i},t_{i_{1}}^{-1}]][w(y_1,\ldots ,y_r),x^{l_i}]}$$
$$ =\overline {\prod_{i=1}^k[w(y_1,\ldots ,y_r),x^{l_i}]}\ \ \ \ \ \ (since\ \
x^{l_i}\in V^{*}(E))$$
$$ = \overline {[w(y_1,\ldots ,y_r),\prod_{i=1}^kx^{l_i}]}$$
$$ = \overline {[w(y_1,\ldots ,y_r),x^n]}\ \ .\ \ \ \Box$$
\newpage
{\bf Lemma 3.5}

Let $B$ be a subgroup of a group $E$ with finite index $n$, and ${\cal V}$ be
the same variety of groups as before. If $x\in V^{*}(E)\cap V(E) $ , then for
all $y_1,\ldots ,y_r\in B$
$$[w(y_1,\ldots ,y_r),x^n]\in V(B)\ \ .$$
In particular, if $L\subseteq V^{*}(E)\cap V(E)$ , then $[W(B),L^n]\subseteq
V(B)\ .$\\
{\bf Proof.}

 Clearly for all $x_i\in E$ ,
$$ ([w(x_1,\ldots ,x_r),x_{r+1},x_{r+2}])\tau =[w(x_1\tau ,\ldots ,x_r\tau
),x_{r+1}\tau ,x_{r+2}\tau ]=1\ \ \ in\ \ B/V(B)\ .$$
So, by Proposition 3.4 for all $x\in V^{*}(E)\cap V(E)$
$$ 1=(x)\tau =[w(y_1,\ldots ,y_r),x^n]\ \ \ in \ \ B/B(V)\ ,$$
which implies that $[w(y_1,\ldots ,y_r),x^n]\in V(B) $ , for all $y_i\in B$ .
This also implies that $[W(B),L^n]\subseteq V(B)$ , for any subgroup $L$ of
$V^{*}(E)\cap V(E)\ .\ \ \Box$

 Now we are able to state and prove the main result of this section.\\
{\bf Theorem 3.6}

 Let ${\cal V}$ be the variety of groups defined by the commutator word
$[w,x,y]$ as before, and $G$ be a finite group with a ${\cal V}$-stem cover
$$ 1\longrightarrow L\longrightarrow G^{*}\longrightarrow G\longrightarrow 1
\ \ .$$
If $H$ is a subgroup of index $n$ in $G$ and $H\cong B/L$ , for some subgroup
$B$ of $G^{*}$ . Then $[W(B),{\cal V}M(G)^n]$ is isomorphic to a subgroup of
${\cal V}M(H)\ .$\\
\newpage
{\bf Proof.}

By the definition of ${\cal V}$-stem cover we have
$$ L\cong {\cal V}M(G)\ \ \ \ \ {\it and}\ \ \ \ L\leq V(G^{*})\cap
V^{*}(G^{*})\ .$$ Clearly $[G^{*}:B]=n$ , and so by Proposition 3.4,
$$[W(B),L^n]\subseteq V(B)\ .$$
Clearly $L$ is marginal in $B$ and since the vatiety ${\cal V}$ has the
Schur-Baer property, Corollary 3.2 implies that $L\cap V(B)$ is isomorphic
to a subgroup of ${\cal V}M(H)$ . We also have
$$ [W(B),{\cal V}M(G)^n]\subseteq L\cap V(B)\ .$$
Hence $[W(B),{\cal V}M(G)^n]$ is isomorphic to a subgroup of ${\cal V}M(H)$ ,
which completes the proof.\ \ $\Box$

 Now the above theorem has the following important corollaries, which
generalize Jones-Wiegold's Theorem [ 3,4 ] and Schur's Theorem [4,9].\\
{\bf Corollary 3.7}

 Let $G$ be a finite group and $G_p$ be a Sylow p-subgroup of $G$ . If $G$ has
a ${\cal V}$-stem cover  $1\rightarrow L\rightarrow G^*\rightarrow
G\rightarrow 1\ , $ such that $G_p\cong B/L\ .$ Then $[W(B),{\cal V}M(G)_p]$ is
isomorphic to a subgroup of ${\cal V}M(G_p)$ .\\
{\bf Proof.}

 Let the group $G$ be of order $p^{\alpha}n$ , where
$(p,n)=1\ ,$ then $[G:G_p]=n\ .$ Using Lemma 2.2 , ${\cal V}M(G)$ has an order
dividing a power of $|G|=p^{\alpha}n\ ,$ and so we have ${\cal V}M(G)^n=
{\cal V}M(G)_p\ .$ Now the result follows from Theorem 3.6 .\ \ $\Box$\\
{\bf Corollary 3.8}

 By the assumptions and the notations of Theorem 3.6, let ${\cal V}M(G)^{[n]}$
be the set of all elements of ${\cal V}M(G)$ of order coprime to $n$ , then
$[W(B), {\cal V}M(G)^{[n]}]$ is isomorphic to a subgroup of ${\cal V}M(H)\ .$\\
{\bf Proof.}

 Clearly, if $G$ is a group then $G^{[n]}$ is a subgroup of $G^n$ . Now the
result holds using Theorem 3.6 .\ \ $\Box$ \\
{\bf Corollary 3.9}

 By the assumptions and the notations of Theorem 3.6 ,let $d$ and $e$ be the
exponents of $[W(B),{\cal V}M(G)]$ and ${\cal V}M(H)$ , respectively. Then $d$
divides $ne$ . In particular, $d$ divides $n$ whenever $H$ is cyclic.\\
{\bf Proof.}

 By the assumption, we have $${\cal V}M(G)\cong L\leq
V^*(G^*),\ \ L\leq B,\ \ {\it and}\ \ {\cal V}M(G)\leq V^*(B)\ .$$
So by Turner-Smith's Theorem [10], it follows that
$$ [W(B),{\cal V}M(G)^n]=[W(B),{\cal V}M(G)]^n\ .$$
Therefore, by Theorem 3.6 ,
$ [W(B),{\cal V}M(G)]^{ne}\subseteq {\cal V}M(H)^e=1\ .$\\
Hence $d$ divides $ne$ .\ \ $\Box$\\

\begin{center}{\bf 4. S{\small\bf CHUR'S} T{\small\bf YPE} T{\small\bf HEOREM}}
\end{center}

 In 1904 , I.Schur [9] proved that if $G$ is a finite group and $G_p$ is a
Sylow $p$-subgroup of $G$ then the Sylow $p$-subgroup of the
Schur-multiplier of $G$, i.e. $M(G)_p$ , is isomorphic to a subgroup of the
Schur-multiplier of its Sylow $p$-subgroup, $M(G_p)$ .

 In this section by using our result together with a theorem of Moghaddam [6]
(see Theorem 2.3), we prove a theorem of Schur's type with respect to the
vatiety of nilpotent groups of class at most $c\ (c\geq 1)$ , ${\cal N}_c$ .

 Now assumming the notations of Theorem 2.3, we first prove the following
technical lemmas.\\
{\bf Lemma 4.1}

 Let the groups $G_1$ and $G_2$ in Theorem 2.3 be the cyclic groups of
orders $q^m$ and $p^n$ , respectively, where $p,q$ are any prime numbers. Then
$T(Z_{p^n},Z_{q^m})_{c+1}=1\ ,$ when $c+1$ is a prime number, and $\hat
{T}(Z_{p^n},Z_{q^m})_4=1\ .$\\
{\bf Proof.}

 Clearly the tensor product $Z_{p^n}\otimes Z_{q^m}=Z_{(p^n,q^m)}=1$ , by the
definition as in Theorem 2.3, it follows that $\hat
{T}(Z_{p^n}\otimes Z_{q^m})=1\ .$

 Now, if $c+1$ is a prime number then $T(Z_{p^n},Z_{q^m})_{c+1}$ is the direct
sum of the following tensor products:
$$ Z_{p^n}\otimes Z_{q^m}\otimes Z_{t_1}\otimes \cdots \otimes Z_{t_{c-1}}\ \
,$$
where $t_i\in \{p^n,q^m\}\ .$ Clearly $Z_{p^n}\otimes Z_{q^m}=1$ and
so each direct summand is trivial. Hence $T(Z_{p^n},Z_{q^m})=1\ .$\\
{\bf Lemma 4.2}

 Let $G=H_1\times H_2$ be the direct product of two finite $p$ and $q$-groups
$H_1$ and $H_2$, respectively. Then ${\cal N}_cM(G)_p\cong {\cal N}_cM(G_p)\
,$ when $c+1$ is any prime number or $c+1=4$ .\\
{\bf Proof.}

 Clearly ${\cal N}_cM(G_p)={\cal N}M(H_1)\ .$ By Theorem 2.3, we have
$$ {\cal N}_cM(H_1\times H_2)\cong {\cal N}_cM(H_1)\oplus {\cal
N}_cM(H_2)\oplus T(H_2,H_1)_{c+1}\ ,$$
where $c+1$ is a prime number. If $c+1=4$ , then the last term was denoted by
$\hat {T}(H_2,H_1)_4\ .$ Clearly any finite abelian p\_group is isomorphic to
some direct sum of cyclic p-groups $Z_{p^{t_i}}$ . Hence using Lemma 4.1 and
the definitions of $T(H_2,H_1)_{c+1}$ and $\hat {T}(H_2,H_1)_4$ , it follows
that both of these are trivial, which gives
$$ {\cal N}_cM(H_1\times H_2)\cong {\cal N}_cM(H_1)\oplus {\cal N}_cM(H_2)\ \
,$$
where $c+1$ is $4$ or any prime number. Now, by Lemma 2.2, ${\cal N}_cM(H_i)$
is a $p$-group and a $q$-group for $i=1,2$ , respectively. Thus
$${\cal N}_cM(G)_p\cong {\cal N}_cM(H_1)={\cal N}_cM(G_p)\ \ .\ \ \Box$$

 Now we state and prove the main result of this section, which is of Schur's
type.\\
{\bf Theorem 4.3}

 Let ${\cal N}_c$ be the variety of nilpotent groups of class at most $c$ and
$G$ be a finite nilpotent group with a Sylow $p$-subgroup $G_p$ , say. If
$c+1$ is any prime number or 4, then
$$ {\cal N}_cM(G)_p\cong {\cal N}_cM(G_p)\ \ .$$
{\bf Proof.}

 It is known that a finite nilpotent group is the direct product of it's Sylow
$p$-subgroups. So
$$ G=G_{p_{1}}\times G_{p_{2}}\times \cdots \times G_{p_{r}}\ \ .$$
Now the result follows by induction on $r$ and applying Lemma 4.2 .\ \ $\Box$

 The following corollary gives a criterion for a finite nilpotent group to
have trivial Baer-invariant with respect to the variety of nilpotent groups of
class at most $c\ (c\geq 1)$.\\
{\bf Corollary 4.4}

 Let ${\cal N}_c$ be the variety of nilpotent groups of class at most $c$ and
$G$ be a finite nilpotent group. If every Sylow $p$-subgroup of $G$ has
trivial Baer-invariant, then so has $G$ .\\
{\bf Proof.}

 Clearly the variety ${\cal N}_c$ has Schur-Baer property, so by Lemma 2.2,
${\cal N}_cM(G)$ has order dividing a power of $|G|$ . Now  Theorem 4.3 gives
the result.\\

\begin{center}  {\bf 5. A C{\small\bf OUNTEREXAMPLE }}  \end{center}

 In this final section by giving an example we show that the result obtained
in the previous sections can not be further generalized. Indeed, take
the symmetric group
$$ S_3=<x,y\ |\ x^2,y^3 ,xyxy>\ \ ,$$
which is isomorphic to the semidirect product of two cyclic groups of orders
two and three. So we have the following free presentation
$$ 1\lra R=<x^2,y^3,xyxy>^F\lra F=<x,y>\lra S_3\lra 1 \ \ .\ \ \ \ \ (*)$$
It is known, by [4] , that $M(S_3)=1$ , so by the definition of
Schur-multiplier we have
$$ R\cap \ga_2(F)=[R,F] $$
Clearly the Bear-invariant of $S_3$ , with respect to the variety of nilpotent
groups of class 2 , is as follows:
$$ {\cal N}_2M(S_3)\cong \frac{R\cap \ga_3(F)}{[R,F,F]}\ \ \ .$$
Now we shall prove that one of the following embeddings does not hold
$$ ({\cal N}_2M(S_3))_2\hookrightarrow {\cal N}_2M((S_3)_2)\ \ \ \ \ \ or\ \ \
\ \ \ ({\cal N}_2M(S_3))_3\hookrightarrow {\cal N}_2M((S_3)_3)\ \ \ .$$
clearly ${\cal N}_2M((S_3)_2)$ and ${\cal N}_2M((S_3)_3)$ are both
trivial, so it is enough to show that either ${\cal N}_2M(S_3)_2$ or ${\cal
N}_2M(S_3)_3$ is non-trivial. But, as ${\cal N}_2$ is a Schur-Baer variety,
Lemma 2.2 implies that $|{\cal N}_2M(S_3)|$ is a 6-number. Hence we only need
to show that ${\cal N}_2M(S_3)$ is non-trivial, which we deal with in the rest
of the section.

 The following technical lemmas shorten the proof of the main theorem. We keep
the above assumption and notation throughout the rest of the paper.

 The proof of the following lemmas is a routine commutator calculation.\\
\newpage
{\bf Lemma 5.1}

 $(i)\ [y,x,y]\in R$ , and $[y,x,y],[y,x,x,y]\in R\cap \ga_2(F)=[R,F]$ ,

 $(ii)\ [y,x,y,a_1,\ldots ,a_n]\in [R,F,F]$ , for all $a_i\in F$ and $n\in {\bf
N}$ ,

 $(iii)\ [y,x,x,y,b_1,\ldots ,b_m]\in [R,F,F]$ , for all $b_i\in F$ and $m\in
{\bf N}$ .\\
{\bf Lemma 5.2}

 The following identities hold, modulo
$$ T_1=<[y,x,y,a_1\ldots ,a_n],[y,x,x,y,b_1,\ldots ,b_m]\ |\ a_i,b_j\in
F,n,m\in {\bf N}>:$$
$(i)\ [x^2,y,y]\equ [y,x,y]^{-2}[x,y,x,y]$\\
$(ii)\ [x^2,y,x]\equ [x,y,x]^{[x,y]}[x,y,x,x]^{[x,y]}[x,y,x]$\\
$(iii)\ [y^3,x,y]\equ[y,x,y]^3$\\
$(iv)\ [y^3,x,x]\equ [y,x,x]^{[y,x]^2}[y,x,x]^{[y,x]}[y,x,x]$\\
$(v)\ [xyxy,x,y]\equ [y,x,y]^2[y,x,x,y]$\\
$(vi)\ [xyxy,x,x]\equ [y,x,x]^{[y,x]}[y,x,x,x]^{[y,x]}[y,x,x]$\\
$(vii)\ [xyxy,y,y]\equ [y,x,y]^{-2}[x,y,x,y]$\\
$(viii)\ [xyxy,y,x]\equ [x,y,x]^{[x,y]}[x,y,x,x]^{[x,y]}[x,y,x]$\\
{\bf Lemma 5.3}
$$ [y,x,y]\ ,\ [x,y,x,y]\ ,\ [y,x,x,y]\in [R,F,F]\ \ . $$
{\bf Proof.}

 By Hall-Witt's identities we have
$$ [[y,x],y,x]^{y^{-1}}[y^{-1},x^{-1},[y,x]]^x[x,[y,x]^{-1},y^{-1}]^{[y,x]}=1\
\ . \ \ \ \ \ \ (1) $$
By Lemma 5.1 (ii),
$$ [[y,x],y,x]^{y^{-1}}\in [R,F,F]\ \ .$$
We also have,
$$ [y^{-1},x^{-1},[y,x]]=[y^{-3}y^2,x^{-2}x,[y,x]]\equ [y^2,x,[y,x]]\ \ \ \ \ \
(mod\ [R,F,F])$$
$$\ \ \ \ \ \ \ \ \ \ \equ [y,x,[y,x]]^2\equ 1\ \ \ \ \ \ (mod\ [R,F,F])\ \ .$$
Therefore, by (1)
$$ [y^{-1},x^{-1},[y,x]]^x\in [R,F,F]\ \ \ \Longrightarrow \ \ \
[x,[y,x]^{-1},y^{-1}]\in [R,F,F]\ \ .$$
But
$$ [x,[y,x]^{-1},y^{-1}]\equ [[x,y,x]^{-1},y^{-1}]\equ
[[x,y,x]^{-1},y^{-3}y^2]$$ $$ \equ [x,y,x,y^2]^{-1}\equ [x,y,x,y]^{-2}\ \ \ \ \
\ (mod\ [R,F,F])\ ,$$
hence
$$ [x,y,x,y]^2\in [R,F,F]\ \ .$$
Using Lemma 5.2 (i) , (ii) , $\ \ [x,y,x,y]^3\in [R,F,F]$ . Thus
$$ [x,y,x,y]\in [R,F,F]\ \ \ .$$
Therefore by Lemma 5.2 (i) , (ii) , (v) , we have
$$ [y,x,y]\ \ \ {\it and} \ \ \ [y,x,x,y]\in [R,F,F]\ \ .\ \  \Box$$

 The following lemma can be proved by using Lemmas 5.1,5.2, and 5.3 :\\
\newpage
{\bf Lemma 5.4}

 Let
$$ T_2=<T_1\ ,\ [y,x,x,[f_1,f_2,f_3],a_1,\ldots ,a_n]\ ,\
[y,x,x,[f'_1,f'_2],b_1,\ldots ,b_m]\ ,\ [y,x,y]\ ,$$ $$\ [y,x,x,y]\ ,\
[x,y,x,y]\ |\ f_1,f_2,f_3,f'_1,f'_2,a_i,b_j\in F\ ,\ n,m\in {\bf N}>\ .$$
Then the following identities hold, modulo $T_2$ :\\
$(i)\ [x^2,y,y]\equ [y^3,x,y]\equ [xyxy,x,y]\equ [xyxy,y,y]\equ 1 $\\
$(ii)\ [x^2,y,x]\equ
[y,x,x]^{-2}[y,x,x,x]^{-1}[y,x,x,[y,x]]^3[y,x,x,x,[y,x]]\equ [xyxy,y,x] $\\
$(iii)\ [y^3,x,x]\equ [y,x,x]^3[y,x,x,[y,x]]^3$\\
$(iv)\ [xyxy,x,x]\equ [y,x,x]^2[y,x,x,x][y,x,x,[y,x]][y,x,x,x,[y,x]]$\\
{\bf Lemma 5.5}

 Let
$$ T_3=<T_2\ ,\ [y,x,x]^3\ ,\ [y,x,x,[y,x]]^3\ ,\
[y,x,x,[y,x]]^2[y,x,x,x,[y,x]]\ ,$$
$$ [y,x,x]^2[y,x,x,x][y,x,x,[y,x]]^{-1}\ ,\ [y,x,x,a_1,\ldots ,a_r]^3\ ,$$
$$ [y,x,x,b_1,\ldots ,b_s]^2[y,x,x,x,c_1,\ldots ,c_t]\ ,\
[y,x,x,x,[y,x],d_1,\ldots ,d_u]\ ,$$
$$ |\ a_i,b_j,c_k,d_l\in F\ ,\ r,s,t,u\in {\bf N}>\ \ .$$
Then $T_3^F=T_3$ and $[R,F,F]\leq T_3$ .\\
{\bf Proof.}

 By the definitions of $T_2$ and $T_3$ , one can easily show that $\ T_3^F=T_3$
. By Lemma 5.4 , every generator of $[R,F,F]$ is in $T_3$ , hence
$$ [R,F,F]\leq T_3 \ \ .\ \ \Box$$

 Now we are in a position to prove the main theorem of this section.\\
{\bf Theorem 5.6}
$$ {\cal N}_2M(S_3)\neq 1 \ \ \ .$$
{\bf Proof.}

 Clearly the identity map $id:R\cap \ga_3(F)\lra R\cap \ga_3(F) $ is a
homomorphism, and by Lemma 5.5, $id([R,F,F])\leq T_3$ . Also by $(*)$, we have
$$ T_3\leq R\cap \ga_3(F)\ \ .$$
It follows that the induced map
$$ id^*:\frac{R\cap \ga_3(F)}{[R,F,F]}\lra \frac{R\cap \ga_3(F)}{T_3} $$
is also a homomorphism. By Lemma 5.5, $[y,x,x,[y,x]]\not\in T_3$ , but
$[y,x,x,[y,x]]\in R\cap \ga_3(F)$ so
$$ id^*(\ol{[y,x,x,[y,x]]})\neq 1 $$
Thus
$$ 1\neq \ol{[y,x,x,[y,x]]}\in \frac{R\cap \ga_3(F)}{[R,F,F]}\cong {\cal
N}_2M(S_3)\ \ \ .$$
Hence
$$ {\cal N}_2M(S_3)\neq 1\ \ .\ \ \ \ \Box$$

 Now the following corollary gives the assertion of our counterexample.\\
{\bf Corollary 5.7}
$$ ({\cal N}_2M(S_3))_2\not \hookrightarrow {\cal N}_2M((S_3)_2) $$
or
$$ ({\cal N}_2M(S_3))_3\not \hookrightarrow {\cal N}_2M((S_3)_3)\ \ \ .$$
{\bf Proof.}

 Clearly ${\cal N}_2M((S_3)_2)=1={\cal N}_2M((S_3)_3)$ . By Theorem 5.6
$$ ({\cal N}_2M(S_3))_2\neq 1\ \ \ \ or\ \ \ \ ({\cal N}_2M(S_3)_3)\neq 1\ \
.$$
Hence the result holds.\ \ $\Box$\\
{\bf Remark}

 Clearly $S_3$ is a non-nilpotent solvable group, so the above corollary shows
that the results of section 4 can not be held for solvable groups, in general.
Also it shows that the Schur's Theorem [4,9] and Jones-Wiegold's
Theorem [3,4] can not be held for solvable groups with respect to the
variety of nilpotent groups of class $c\ ,\ c\geq 2$ .

\end{document}